\documentclass[smallextended,numbook,envcountsame]{svjour3}
\usepackage{graphicx,amsmath,amssymb,theorem,mathptmx}
\usepackage[usenames]{color}
\smartqed
\numberwithin{equation}{section}

\newcommand{\Real}{\mathbb R}

\journalname{Mathematische Annalen}
\begin{document}

\title{Universal bounds and semiclassical estimates for eigenvalues of 
abstract Schr\"odinger operators}%
\author{Evans M. Harrell II \and Joachim Stubbe}

\institute{E. M. Harrell \at School of Mathematics, 
Georgia Institute of Technology, Atlanta GA 30332-0160, USA\\
\email{harrell@math.gatech.edu}\\
\and{J. Stubbe \at 
EPFL, IMB-FSB, Station 8, CH-1015 Lausanne, Switzerland\\
\email{Joachim.Stubbe@epfl.ch}
}}

\date{Received: date / Accepted: date}
\maketitle

\begin{abstract}
We prove trace inequalities for a self-adjoint operator on an
abstract Hilbert space. These inequalities lead to universal
bounds on spectral gaps and on moments of eigenvalues
$\{\lambda_k\}$
that are analogous to those known for Schr\"odinger operators 
and the Dirichlet Laplacian, on which the operators of interest are 
modeled.  In addition we produce inequalities that are new
even in the model case.  These include a family of differential 
inequalities for generalized Riesz means and theorems stating that
arithmetic means of $\{\lambda_k^p\}_{k=1}^n$ for $p \le 3$ are universally 
bounded from above by 
multiples of the 
geometric means, $\left(\prod_{k=1}^n{\lambda_k}\right)^{1/n}$.
For
Schr\"odinger operators and the Dirichlet Laplacian these bounds
are Weyl-sharp, i.e., saturated by the standard semiclassical estimates
as $n \rightarrow \infty$.


\end{abstract}
\section{Introduction}

Universal spectral bounds for Laplace and Schr\"odinger operators,
i.e., bounds that control eigenvalues with expressions that do not
depend on the specific geometry of the domain or on details of the
potential (cf. \cite{As}), can be derived from fundamental
identities involving traces of operators and their commutators
\cite{HaSt}. This insight has proved useful both for unifying
numerous previously known inequalities of this kind and for
discovering new ones \cite{HaSt,AsHe,HaHe,HaHe2,LePa}. Related
methods have also been used to obtain control on the spectrum of
Laplace and Schr\"odinger operators in terms of curvature
\cite{CY1,CY2,Ha07,EHI}. 
Many of the universal inequalities related to trace identities
are sharp in the sense that they are saturated for particular
examples:  For the Schr\"odinger operators treated in \cite{HaSt}
the upper bounds on eigenvalue gaps $\lambda_{n+1} - \lambda_n$
become identities for all $n$ in the case of the harmonic
oscillator, while for Laplacians on embedded manifolds discussed
in \cite{Ha07,EHI} all of the gap bounds become identities for
embedded spheres. More recently, in some circumstances (e.g.,
\cite{Lap,LW1,LW2,HaHe2}) universal bounds on moments of
eigenvalues have been connected to ``semiclassical'' theorems
about the spectrum such as asymptotic behavior as the index $k
\rightarrow \infty$ and nonasymptotic bounds in the spirit of the
Berezin-Li-Yau inequality \cite{Be,LiYa}.

 One of the motivations of this work is to
sharpen the understanding of moments of eigenvalues and of Riesz
means of the spectrum, which plays something like the role of a
dual version of moments. Among the applications of our analysis
will be a family of differential inequalities for functions
determined by the spectrum, extending the analysis of
\cite{HaSt,HaHe}. By Legendre duality as in \cite{HaHe} these
imply bounds on ratios of averages of eigenvalues.  We also
introduce a novel type of inequality relating arithmetic and
geometric means of eigenvalues.

A second motivation is to better unify the subject of universal
bounds with formally analogous semiclassical spectral theorems.

 In the next section we present some more
abstract versions of the essential trace identity of \cite{HaSt}
for a class of self-adjoint operators $H$ enjoying algebraic
properties modeled on those of of Schr\"odinger operators.  
We also identify a special
family of functions for which $tr(H)$ can be sharply
controlled, {\it viz.}:
\begin{definition}\label{SJ} 
 Let $H$ be a self-adjoint operator and let $J
\subset \sigma(H)$ be a distinguished subset of the spectrum.  We
let ${\hat J}$ denote the smallest closed interval containing $J$.
A $C^1$ function $f: {\hat J} \rightarrow \mathbb{R}$ belongs to
the set ${\mathfrak S}_J$ of {\it trace-controllable functions}
provided that on ${\hat J}$,
\begin{description}[H1.~]
\item[H1.] $f(\lambda) \ge 0$; \item[H2.] $f^\prime(\lambda) \le
0$; \item[H3.] $f^\prime(\lambda)$ is concave; \item[H4.] If
$\sup(J) < \infty$, there exists $a > \sup(J)$ such that
$$g_f(x) := (a - \lambda)^3 \frac{d}{d\lambda} \left( \frac{f(\lambda)}{(a -
\lambda)^2}\right) = 2 f(\lambda) + f^\prime(\lambda) (a -
\lambda)$$ is nondecreasing in $\lambda$; \item[H5.] $tr(P_J(H)
f(H)) < \infty$, where $P$ denotes the spectral projector for the
set $J$.
\end{description}
\end{definition}

For reasons of parsimony we shall sometimes assume only a subset
of the hypotheses in the statements of some theorems.

The model situation is that at least the lower part of the
spectrum consists of eigenvalues $\lambda_1 < \lambda_2 \le \cdots$
and $J = \{\lambda_1, \dots, \lambda_n\}$. Among familiar
functions in ${\mathfrak S}_J$ we mention $\exp(-t \lambda)$ and
$(z - \lambda)^p$ with $p \ge 2$. Note that conditions H1--H4 are
preserved by multiplication $f(\lambda),g(\lambda) \rightarrow
f(\lambda)g(\lambda)$ and that conditions H1--H3 are preserved by
compositions in the form $f(\lambda),g(\lambda) \rightarrow
f(-g(\lambda))$.

The only condition in Definition~\ref{SJ} that may not be familiar
is H4, so we observe some sufficient conditions for its validity,
depending on some elementary facts about concavity, 
in particular,

\begin{proposition}\label{cvxprop}
If the function $h(x)$ is concave for $0<x<x_0$, then
$$
x h(x) - 2 \int_0^x{h(s) ds}
$$
is concave on the same interval.
\end{proposition}

We observe that this is immediate when $h \in C^2$ by a
calculation of the second derivative.  We give a proof without
this hypothesis.
\begin{proof}
Recall that a function f is concave on an interval $I$ iff its right
and left derivatives exist at all interior points of $I$ and $f'(x)$
is nonincreasing, 
in the extended sense that if the right and left derivatives differ at $x$, then $f_r'(x) < f_{\ell}'(x)$.
(For this and other basic facts about concave functions see chapter 5 of
\cite{We}.)
Therefore we compare the derivative
 $\frac{d}{dx} (x h(x) - 2 \int_0^x h(s) ds ) = xh^\prime(x) - h(x)$ at $a$ (right derivative) and $a+\delta$
(left derivative), for $\delta > 0$.  
(We shall not complicate the notation by distinguishing right and left derivatives in the
following calculation.)
\begin{align*}
    &(x+\delta)h^\prime(x+\delta) - h(x+\delta) - (x h^\prime(x) - h(x) )\\
    &\quad = \delta h^ \prime(x+\delta) + x (h^ \prime(x+\delta) - h^\prime(x)) - (h(x+\delta) -
    h(x))\\
     &\quad =x (h^\prime(x+\delta) - h^\prime(x)  + \left(\delta h^ \prime(x+\delta) - ((h(x+\delta) -
     h(x))\right)\\
     &\quad  \leq \delta h^ \prime(x+\delta) - ((h(x+\delta) - h(x)).
\end{align*}
By the mean value theorem of convex functions,
for some $y \in (x,x+\delta)$, $(h_{\ell}^\prime(y) \cdot \delta
\le h(x+\delta) - h(x)) \le (h_{r}^\prime(y) \cdot \delta$, and
thus the final term is $\le 0$.\qed
\end{proof}

\begin{corollary}\label{gfunction}
Suppose that $f$ satisfies hypotheses {\rm H2} and {\rm H3} of 
Definition~{\rm\ref{SJ}} and that

{\rm H4$'$}  There exists $a > \sup(J)$ such that
$$g_f'(\sup(J) = f^\prime(\sup(J)) + f^{\prime\prime}(\sup(J)) (a - \sup(J)) \ge 0.$$
Then $f$ satisfies hypothesis {\rm H4}.
\end{corollary}

\begin{proof}
We study the function $g_f$ occurring in Hypothesis H4. With the
change of variable $a-\lambda \rightarrow x$, $h(x) :=
f^\prime(\lambda)$ satisfies the conditions of 
Proposition~\ref{cvxprop}, so $x h(x) - 2 \int_0^x{h(s) ds}$ is concave for
positive $x$.  But this expression evaluates to $g_f(\lambda) -
g_f(a)$, establishing that $g_f$ is concave for $\lambda \le a$.
Therefore $g_f^\prime$ is nonincreasing.  At the same time we know
by H4$'$ that $g_f'(\sup(J) \ge 0$, so it follows that
$g_f'(\lambda) \ge 0$ on ${\hat J}$.\qed
\end{proof}

\section{Abstract trace inequalities}

We consider a self-adjoint operator $H$ with domain
$\mathcal{D}_H$ on a Hilbert space $\mathcal{H}$ with scalar
product $\langle\cdot,\cdot\rangle$. We suppose that
$H$ has nonempty point spectrum, and that $\mathcal{J}$ is a
finite-dimensional subspace of $\mathcal{H}$ spanned by an
orthonormal set $\left\{\phi_j\right\}$ of eigenfunctions of $H$.
We further let $P_{A}$ denote the spectral projector associated
with $H$ and a Borel set $A$, and $J := \left\{\lambda_j: H \phi_j
= \lambda_j \phi_j\right\}$. We refer to \cite{ReSi} for
terminology, notation, and details about the spectral theorem.

\begin{theorem}\label{main theorem 1} Let $H$ and $G$ be
self-adjoint operators with domains $\mathcal{D}_H$ and
$\mathcal{D}_G$ such that $G(\mathcal{J})\subseteq
\mathcal{D}_H\subseteq \mathcal{D}_G$.
Then, for any real-valued $C^1$-function $f$ the derivative $f'$
of which is a concave function $($i.e., Hypothesis~{\rm H3} of 
Definition~{\rm\ref{SJ})},
\begin{align}\label{T0}
    &\frac1{2}\sum_{\lambda_j\in J}f'(\lambda_j)\,\langle[H,G]\phi_j,[H,G]\phi_j\rangle +  f(\lambda_j)\,\langle[G, [H,G],]\phi_j,\phi_j\rangle\\
    &\qquad \leq
    \sum_{\lambda_j\in J} \int{
    \big(f(\lambda_j)+\frac1{2}f'(\lambda_j)(\kappa-\lambda_j)\big)(\kappa-\lambda_j)|\langle G\phi_j,dP_{\kappa}
    P_{{\mathcal J}^c} G\phi_j \rangle|^2}\nonumber
\end{align}
In case the spectrum of $H$ is purely discrete, we may write the
inequality as
\begin{align}\label{T1}
    &\frac1{2}\sum_{\lambda_j\in J}f'(\lambda_j)\,\langle[H,G]\phi_j,[H,G]\phi_j\rangle +  f(\lambda_j)\,\langle[G, [H,G],]\phi_j,\phi_j\rangle\\
    &\qquad \leq
    \sum_{\lambda_j\in J}\sum_{\lambda_k\in J^c}
    \big(f(\lambda_j)+\frac1{2}f'(\lambda_j)(\lambda_k-\lambda_j)\big)(\lambda_k-\lambda_j)|\langle G\phi_j,\phi_k\rangle|^2.\nonumber
\end{align}
\end{theorem}

The proof will use the following lemma:

\begin{lemma}\label{cvxlem} Let $f\in C^1(\Real)$ such that $f'$ is a concave
function. Then for all $x,y\in \Real$
\begin{equation}\label{f}
    \frac{f(y)-f(x)}{y-x}\geq \frac1{2}f'(y)+\frac1{2}f'(x).
\end{equation}
\end{lemma}

\begin{proof}
By the fundamental theorem of calculus and the concavity of $f'$
we have
\begin{align*}
    \frac{f(y)-f(x)}{y-x}&=\int_{0}^{1}f'((1-t)x+ty)\,dt\\&\geq
    \int_{0}^{1}(1-t)f'(x)+tf'(y)\,dt
    =\frac1{2}f'(y)+\frac1{2}f'(x).\qquad \qed
\end{align*}
\end{proof}

\begin{proof}
We begin with an observation that is an abstract version of what
is known in quantum theory as the oscillator-strength sum rule of
Thomas, Reiche, and Kuhn \cite{BeJa}:  By a straightforward
calculation, the self-adjoint operators $(H,G)$ satisfy
\begin{equation}\label{TRK}
\langle[G, [H,G]]\phi_j,\phi_j\rangle = 2 \langle (H - \lambda_j)
G\phi_j, G\phi_j\rangle,
\end{equation}
which, with the spectral resolution, equals $2 \int (\kappa  -
\lambda_j)\langle dP_{\kappa} G\phi_j, G\phi_j\rangle.$ Thus
\begin{equation}\label{sumrule 1}
   \frac1{2}\langle[G,[H,G]]\phi_j,\phi_j\rangle =
   \int_{}{(\kappa-\lambda_j)dG_{j \kappa}^2},
\end{equation}
where
$dG_{j \kappa}^2:=|\langle G\phi_j,dP_{\kappa} G\phi_j\rangle|$.
When $\kappa = \lambda_k$ for $\phi_k \in {\mathcal J}$, we also
write the discrete matrix elements as $G_{j k}:=\langle
G\phi_j,\phi_k\rangle$.

Multiplying by $f(\lambda_j)$ and summing over $\lambda_j$  in
$J$, we get
\begin{align}\label{1a}
    &\frac1{2}\sum_{\lambda_j\in J}f(\lambda_j)\,\langle[G, [H,G]]\phi_j,\phi_j\rangle 
     =\sum_{\lambda_j\in J}\int_{}
    f(\lambda_j)(\kappa-\lambda_j)dG_{jd}^2\\
    &\qquad =\sum_{\lambda_j\in J}\sum_{\lambda_k\in J}
    f(\lambda_j)(\lambda_k-\lambda_j)G_{jk}^2+\sum_{\lambda_j\in J}\int_{\kappa \in J^c}
    f(\lambda_j)(\kappa-\lambda_j)dG_{j \kappa}^2.\nonumber
\end{align}
Using the symmetry of the matrix elements $G_{jk}$ we rewrite the
first double sum as follows:
\begin{align*}
    \sum_{\lambda_j\in J}\sum_{\lambda_k\in J}
    f(\lambda_j)(\lambda_k-\lambda_j)G_{jk}^2
    &=\frac1{2}\sum_{\lambda_j\in J}\sum_{\lambda_k\in J}
    (f(\lambda_j)-f(\lambda_k))(\lambda_k-\lambda_j)G_{jk}^2\\
    &=-\frac1{2}\sum_{\lambda_j\in J}\sum_{\lambda_k\in J}
    \frac{f(\lambda_k)-f(\lambda_j)}{\lambda_k-\lambda_j}(\lambda_k-\lambda_j)^2G_{jk}^2.
\end{align*}
Applying Lemma~\ref{cvxlem} and once again using the symmetry of
the $G_{jk}$ we get
\begin{equation}\label{1b}
     \sum_{\lambda_j\in J}\sum_{\lambda_k\in J}
    f(\lambda_j)(\lambda_k-\lambda_j)G_{jk}^2
    \leq
     -\frac1{2}\sum_{\lambda_j\in J}\sum_{\lambda_k\in J}
    f'(\lambda_j)(\lambda_k-\lambda_j)^2G_{jk}^2.
\end{equation}
At the same time, the pair $(H,G)$ satisfies the trace formula
\begin{equation*}
    \int{(\kappa-\lambda_j)^2 dG_{j \kappa}^2}=\langle[H,G]\phi_j,[H,G]\phi_j\rangle.
\end{equation*}
Multiplying by $-\frac1{2}f'(\lambda_j)$ and summing over
$\lambda_j\in J$ we get
\begin{align}\label{1c}
    &-\frac1{2}\sum_{\lambda_j\in J}f'(\lambda_j)\,\langle[H,G]\phi_j,[H,G]\phi_j\rangle \\
    &\quad =-\frac1{2}\sum_{\lambda_j\in J}\int{
    f'(\lambda_j)(\kappa-\lambda_j)^2 dG_{j \kappa}^2}\nonumber\\
    &\quad =-\frac1{2}\sum_{\lambda_j\in J}\sum_{\lambda_k\in J}
    f'(\lambda_j)(\lambda_k-\lambda_j)^2G_{jk}^2\nonumber\\
    &\qquad-\frac1{2}\sum_{\lambda_j\in J}\int_{\kappa \in J^c}
    f'(\lambda_j)(\kappa-\lambda_j)^2 dG_{j \kappa}^2.\nonumber
\end{align}
Combining \eqref{1a}, \eqref{1b} and \eqref{1c} yields the
statement of the theorem.\qed
\end{proof}

\begin{corollary}\label{quadratic function}
If $f(\lambda)=a\lambda^2+b\lambda +c$, then \eqref{T1} holds with
equality. In particular, for any $z$ we have
\begin{align}\label{C1}
    & \sum_{\lambda_j\in J}{(z-\lambda_j)^2\,\langle[G, [H,G]]\phi_j,\phi_j\rangle -
 2(z-\lambda_j)\,\langle[H,G]\phi_j,[H,G]\phi_j\rangle}
    \\
    &\qquad =2 \sum_{\lambda_j\in J}\int_{\kappa \in J^c}
    (z-\lambda_j)(z-\kappa)(\kappa-\lambda_j) d G_{j \kappa}^2.\nonumber
\end{align}
\end{corollary}

\begin{remark}
Equation~\eqref{C1} is a particular case of a more general trace
identity that will be explored in a future work.
\end{remark}

We now add for the first time the assumption that the lower part
of the spectrum of $H$ is discrete and denote this set
$J=\{\lambda_1,\ldots\lambda_n\}$. The following theorem captures
a universal relationship between the lower part of the spectrum
for $j = 1, 2, \dots, n$ and the values of $\lambda_n$ and
$\lambda_{n+1}$.

\begin{theorem} \label{inequalities for bottom of spectrum} Let $H$ and $G$ be
self-adjoint operators with domains $\mathcal{D}_H$ and
$\mathcal{D}_G$ such that $G(\mathcal{J})\subseteq
\mathcal{D}_H\subseteq \mathcal{D}_G$.  Let the subset
$J=\{\lambda_1,\ldots\lambda_n\}$ lie below the rest of the
spectrum of $H$. Then for any
$f \in {\mathfrak S}_J$,  with
$f'(\lambda_n)+f''(\lambda_n)(\lambda_{n+1}-\lambda_n)\geq 0$,
\begin{align}\label{T2}
    &\frac1{2}\sum_{j=1}^{n}\left(f'(\lambda_j)\,\langle[H,G]\phi_j,[H,G]\phi_j\rangle
    + f(\lambda_j)\,\langle[G, [H,G]]\phi_j,\phi_j\rangle\right)\\
    &\qquad \leq
    \frac1{2}\big(f(\lambda_n)+\frac1{2}f'(\lambda_n)(\lambda_{n+1}-\lambda_n)\big)\sum_{j=1}^n\langle[G, [H,G]]\phi_j,\phi_j\rangle.\nonumber
\end{align}
\end{theorem}
\begin{proof}
Consider the right side of \eqref{T1}. Since
$\lambda_k\geq\lambda_{n+1}\geq\lambda_j$ and $f'\leq 0$, we have
\begin{equation*}
    f(\lambda_k)+\frac1{2}f'(\lambda_k)(\lambda_{l}-\lambda_k)\leq
    f(\lambda_k)+\frac1{2}f'(\lambda_k)(\lambda_{n+1}-\lambda_k) = \frac{1}{2} g_f(\lambda_k).
\end{equation*}
To prove \eqref{T2} it suffices to show that $g$ is
nondecreasing so $g_f(\lambda_k)$ can be replaced
with $g_f(\lambda_n)$. As a concave function, $g^\prime$ is nondecreasing
on ${\hat J}$, so this is true because of the assumption
that $g^\prime(\lambda_n)
\ge 0$.\qed
\end{proof}

Under slightly weakened assumptions on $f$ we get the following
(weaker) inequality:
\begin{corollary} \label{weaker inequalities for bottom of spectrum} Let $H$ and $G$ be
self-adjoint operators with domains $\mathcal{D}_H$ and
$\mathcal{D}_G$ such that $G(\mathcal{J})\subseteq
\mathcal{D}_H\subseteq \mathcal{D}_G$.  Let the subset
$J=\{\lambda_1,\ldots\lambda_n\}$ lie below the rest of the
spectrum of $H$. Then, for any real $C^1$-function $f$ satisfying
Hypothesis~{\rm H2} with $f'(\lambda_{n+1})= 0$,
\begin{align} 
    &\frac1{2}\sum_{j=1}^{n}\left(f'(\lambda_j)\,\langle[H,G]\phi_j,[H,G]\phi_j\rangle
    + f(\lambda_j)\,\langle[G, [H,G]]\phi_j,\phi_j\rangle\right)\\
    &\qquad \leq
    \frac1{2}f(\lambda_{n+1})\sum_{j=1}^n\langle[G, [H,G]]\phi_j,\phi_j\rangle.\nonumber
\end{align}
\end{corollary}

In applications to Laplace, Schr\"odinger, and similar
differential operators, the commutators typically simplify as
follows: There are constants $\alpha, \beta,\gamma,$ with
$\beta,\gamma>0$, such that
\begin{equation}\label{commutator assumptions}
    \gamma=\frac1{2}[G,[H,G]],\quad \beta H+\alpha \,\, \ge\, -[H,G]^2.
\end{equation}
(Recall that $- [H, G]^2 = [H,G] [H,G]^* \ge 0$.) For instance,
see \cite{As,AsHe,ChYa,HaSt,Ha07,LePa,Li,Ya}. Therefore we also
may assume without loss of generality that $H$ has only
nonnegative eigenvalues. Indeed, if $\tilde{H}=H+\eta$ for some
real constant $\eta$, then
\begin{equation*}
    \tilde{\lambda}=\lambda+\eta,\quad \tilde{\gamma}=\gamma, \quad
    \tilde{\beta}=\beta, \quad \tilde{\alpha}=\alpha-\beta\eta.
\end{equation*}

For the model case of the Dirichlet Laplacian $H=-\Delta_D$ on a
domain $\Omega$ in $\mathbb{R}^d$ there is a choice of Cartesian
system coordinate system for such that with $G = x_1$, $\alpha =
0, \beta = \frac{4}{d}$, and $\gamma = 1$. In the literature these
same effective constants are often obtained by averaging over all
the coordinates, but the coordinate system can always be chosen to
make this unnecessary.

Values of $\alpha \ne 0$ arise for several reasons. In the case of
a Schr\"odinger operator $H = - \Delta + V({\bf x})$, the
potential energy disappears from all commutators, and the term $-
[H, G]^2$ is typically dominated by the kinetic energy term, i.e.,
$- [H, G]^2 \le \beta (- \Delta)$ rather than $\beta H$.  In an
elementary way, the addition of $\alpha$ can compensate for the
absence of $V$ if, say, the negative part of $V$ is bounded.  Even
if the negative part of the potential $V$ is unbounded, if it lies
in certain function classes, there are constants $a < 1$ and $b <
\infty$ such that for all functions $\varphi$ in the
quadratic-form domain of $H$,
\begin{equation}\label{qfb}
|\left\langle\varphi , V_-({\bf x}) \varphi \right\rangle| \le a
\|\nabla \varphi\|^2 + b  \| \varphi\|^2,
\end{equation}
in which case $ -\Delta + V_- \ge (1-a) (-\Delta) - b$, and
consequently
$$
-\Delta \le \frac{1}{1-a} \left(H + b\right).
$$
Examples of function classes guaranteeing the estimate \eqref{qfb}
are that $V_-$ is a Rollnik potential in three dimensions or that
$V_- \in L^p(\mathbb{R}^d)$ with $p > \frac{d}{2}$ when $d \ge 4$.
For discussion of these conditions refer to \cite[section~X.2.]{ReSiII}.

Another instance where $\alpha \ne 0$ is of interest is in the
case of Laplace or Schr\"odinger operators on hypersurfaces
${\mathcal M}$ in $(\mathbb{R}^d)$.  By letting $G$ be the
Cartesian coordinate $x_1$ in the ambient space, and choosing the
orientation of the coordinate system appropriately (or averaging
over all coordinates),
$$
- [H, G]^2 =  - 4\Delta + h^2(x),
$$
where $h(x)$ is the sum of the principal curvatures at the point
$x \in {\mathcal M}$ and $\Delta$ now denotes the Laplace-Beltrami
operator on ${\mathcal M}$ \cite{Ha07}.  Our estimates therefore
apply to Laplace-Beltrami operators on ${\mathcal M}$ with a
$\alpha = \|h\|_{\infty}^2$.  Schr\"odinger operators on
${\mathcal M}$ will require $\alpha$ to be the sum of this
curvature effect and any contribution owing to the negative part
of $V$.  The situation is analogous for Laplace or Schr\"odinger
operators on manifolds immersed in other symmetric spaces
\cite{EHI}.

Under these conditions Theorem~\ref{inequalities for bottom of
spectrum} is simplified:

\begin{corollary}\label{corollary to inequalities for bottom of spectrum}
If in addition to the assumptions of Theorem~{\rm\ref{inequalities for
bottom of spectrum}} the relations \eqref{commutator assumptions}
hold, then
\begin{equation}\label{C2}
\frac1{n}\sum_{j=1}^{n}\left(f(\lambda_j)+\frac{\beta\lambda_j+\alpha}{2\gamma}f'(\lambda_j)\right)
\leq f(\lambda_n)+\frac1{2}f'(\lambda_n)(\lambda_{n+1}-\lambda_n).
\end{equation}
 If, in addition, the spectrum of $H$ is purely discrete and all
sums over the full spectrum $\sigma(H)$ are finite, then
\begin{equation}\label{C2bis}
\sum_{\lambda_j\in\sigma(H)}\left(f(\lambda_j)+\frac{\beta\lambda_j+\alpha}{2\gamma}f'(\lambda_j)\right)\leq
0.
\end{equation}
\end{corollary}

In the following sections we apply Corollary~\ref{corollary to
inequalities for bottom of spectrum} for appropriate functions
$f(\lambda)$.

\section{Inequalities for moments of eigenvalues}
In this section we prove various inequalities for eigenvalues
under the assumptions of Corollary~\ref{corollary to inequalities
for bottom of spectrum} for appropriate functions $f(\lambda)$,
and we restrict ourselves to operators $H$ with purely discrete
spectrum. As a first result we generalize the result of
\cite{HaSt} on the partition function $\text{tr}(e^{-tH})$.

\begin{proposition}
Let $f$ and $H$ satisfy the assumptions of Corollary~{\rm\ref{corollary to
inequalities for bottom of spectrum}} and suppose that
\begin{equation*}
    \text{tr}\left(f\left(t\left(H+\frac{\alpha}{\beta}\right)\right)\right)
=\sum_{j=1}^{\infty}f\left(t\left(\lambda_j+\frac{\alpha}{\beta}\right)\right)
\end{equation*}
and
\begin{equation*}
   \frac{d}{dt}\; \text{tr}\left(f\left(t\left(H+\frac{\alpha}{\beta}
\right) \right) \right)
=\sum_{j=1}^{\infty}\left(\lambda_j+\frac{\alpha}{\beta}\right)
f\left(t\left(\lambda_j+\frac{\alpha}{\beta}\right)\right)
\end{equation*}
are finite for all $t>0$, then
\begin{equation} 
t \mapsto
t^{\frac{2\gamma}{\beta}}\text{tr}\left(f\left(t\left(H+\frac{\alpha}{\beta}
\right)\right)\right)
\end{equation}
is nonincreasing.
\end{proposition}

\begin{proof}
If $f(\lambda)$ satisfies the assumptions of Theorem~\ref{main
theorem 1}, then $f(t\lambda+ \frac{\alpha}{\beta}))$ satisfies
the same assumptions for any $t>0$. Then inequality \eqref{C2bis}
of Corollary~\ref{corollary to inequalities for bottom of
spectrum} reads as follows:
\begin{equation*}
\text{tr}\left(f\left(t\left(H+\frac{\alpha}{\beta}
\right)\right)\right)+\frac{\beta}{2\gamma}\;
t\frac{d}{dt}\; \text{tr}\left(f\left(t\left(H+\frac{\alpha}{\beta}
\right)\right)\right)\leq 0,
\end{equation*}
which proves the proposition.\qed
\end{proof}

The proposition applies to $f(\lambda)=\lambda^{-p}e^{-\lambda}$,
for any $p\geq 0$ and therefore
\begin{corollary} \label{weigthed heat kernel}
If
\begin{equation*}
    Z_p(t):=\text{tr}
\left( \left(H+\frac{\alpha}{\beta}\right)^{-p}e^{-tH}\right)
=\sum_{j=1}^{\infty}\left(\lambda_j+\frac{\alpha}{\beta}\right)^{-p}e^{-t\lambda_j}
\end{equation*}
is finite for all $t>0$ and $H$ satisfies the assumptions of
Corollary~{\rm\ref{corollary to inequalities for bottom of spectrum}},
then
\begin{equation} 
t \mapsto
Z_p(t)t^{\frac{2\gamma}{\beta}-p}e^{-\frac{\alpha}{\beta}t}
\end{equation}
is nonincreasing.
\end{corollary}

\begin{remark}
In particular, Corollary~\ref{weigthed heat
kernel} shows that
\begin{equation*}
\underset{t\rightarrow 0+}{\lim}Z_p(t)=+\infty
\end{equation*}
for all $p<\frac{2\gamma}{\beta}$.
\end{remark}

As a second application, we shall show that
certain moments of eigenvalues are dominated by their geometric
mean. Let $z>0$ be a parameter to be chosen later and $p>q>0$ such
that $q\leq \min(1,p)$ and $p+q\leq 3$. For $\lambda\in [0,z]$ the
function $f_z(\lambda)$ defined by
\begin{equation}\label{small powers}
    f_z(\lambda):=q\lambda^p-p\lambda^{q}z^{p-q}+(p-q)z^p
\end{equation}
satisfies the assumptions of Theorem~\ref{inequalities for bottom
of spectrum} and Corollary~\ref{corollary to inequalities for
bottom of spectrum} for all $z\in[\lambda_n,\lambda_{n+1}]$.
Indeed, with
\begin{align*}
  f'_z(\lambda) &= pq(\lambda^{p-1}-\lambda^{q-1}z^{p-q}) \\
  f''_z(\lambda) &= pq\lambda^{q-2}\big((p-1)\lambda^{p-q}-(q-1)z^{p-q}\big)\\
 f'''_z(\lambda) &=
 pq\lambda^{q-3}\big((p-1)(p-2)\lambda^{p-q}-(q-1)(q-2)z^{p-q}\big)
\end{align*}
we see that $f''_z(\lambda)\geq 0$ if $q\leq 1$, using the
estimate $(1-q)z^{p-q}\geq (1-q)\lambda^{p-q}$. Furthermore,
$f'''_z(\lambda)\leq 0$ since
$(p-1)(p-2)\lambda^{p-q}-(q-1)(q-2)z^{p-q}\leq
(p-q)(p+q-3)\lambda^{p-q}$. As in the proof of 
Theorem~\ref{inequalities for bottom of spectrum} we show that $g'(z)\geq
0$ if $z\leq \lambda_{n+1}$, and hence
\begin{equation*}
    g(\lambda_k)\leq g(\lambda_n)\leq g(z)=0,
\end{equation*}
provided that $z\in[\lambda_n,\lambda_{n+1}]$. We define
\begin{equation*}
    F_n(z) := \sum_{j=1}^{n}f_z(\lambda_j)+\frac{\beta\lambda_j+\alpha}{2\gamma}f_z'(\lambda_j).
\end{equation*}
Applying Corollary~\ref{corollary to inequalities for bottom of
spectrum} with $0$ as an upper bound, we have $F_n(z)\leq 0$ for
all $z\in[\lambda_n,\lambda_{n+1}]$. Since
\begin{equation*}
    F_n(z)=n(p-q)z^p-npS_n(q)z^{p-q}+nqS_n(p),
\end{equation*}
where
\begin{equation*}
    S_n(r) := \left(1+\frac{\beta
    r}{2\gamma}\right)\frac1{n}\sum_{j=1}^{n}\lambda_j^r+\frac{\alpha r}{2\gamma}\frac1{n}\sum_{j=1}^{n}\lambda_j^{r-1},
\end{equation*}
we see that $F_n(z)$ attains a global nonnegative minimum at
$z=S_n(q)^{\frac{1}{q}}$. Therefore we have the following result:
\begin{theorem}\label{main theorem on moments} Let $p>q>0$ such that
$q\leq \min(1,p)$ and $p+q\leq 3$. Then for all $n$ we have
\begin{equation}\label{means small p}
    S_n(p)^{\frac1{p}}\leq  S_n(q)^{\frac1{q}}.
\end{equation}
In particular, for $0<p\leq 1$ the function $p\mapsto
S_n(p)^{\frac1{p}}$ is nonincreasing, and for all $0<p\leq 3$ we
have
\begin{equation}\label{geometric mean}
    S_n(p)^{\frac1{p}}\leq e^{\frac{\beta}{2\gamma}}G_n
\exp\left(\frac{\alpha}{2\gamma}\frac1{n}\sum_{j=1}^{n}\lambda_j^{-1}\right),
\end{equation}
where
\begin{equation*}
    G_n:=\bigg(\prod_{j=1}^{n}\lambda_j\bigg)^{\frac1{n}}
\end{equation*}
denotes the geometric mean of the first $n$ eigenvalues.
\end{theorem}
Inequality \eqref{geometric mean} is obtained from \eqref{means
small p} by taking the limit $q\rightarrow 0$.
\begin{remark}
We note that $F_{n+1}(\lambda_{n+1})=F_n(\lambda_{n+1})$ for all
$n$.
\end{remark}

\section{Inequalities on generalized Riesz means}

Recently, it was shown in \cite{HaHe} that Riesz means of
eigenvalues of the Laplacian on bounded domains satisfy
differential inequalities, which in turn imply universal
eigenvalue bounds that are sharp in the sense of having the
correct behavior as expected from the Weyl law for $\lambda_n$ as
$n \rightarrow \infty$. We note here how some of the results of
\cite{HaHe} can be extended as a consequence of 
Theorem~\ref{main theorem 1}.

Let $f$ be a function in the class ${\mathfrak S}_J$ such that
$f(1)=f'(1)=0$. We define a generalized Riesz mean by
\begin{equation}\label{generalized Riesz mean}
    R_f(t):=\sum_{j}f(t\lambda_j)\theta(1-t\lambda_j),
\end{equation}
where $\theta(x)=1$ if $x>0$ and zero otherwise and $f$ satisfies
the hypotheses of Theorem~\ref{main theorem 1}.

For simplicity we consider the case $\alpha=0$,
which can always be arranged
by shifting $\lambda_j\rightarrow \lambda_j+\alpha/\beta$.

\begin{corollary}\label{corollary to generalized Riesz means}
If in addition to the assumptions of Theorem~{\rm\ref{inequalities for
bottom of spectrum}} the relations \eqref{commutator assumptions}
hold, then
\begin{equation}\label{R1}
\frac{d}{dt}\;t^{\frac{2\gamma}{\beta}}R_f(t)\leq 0.
\end{equation}
\end{corollary}

In \cite{HaHe} the Riesz means $$R_{\rho}(z) :=
\sum_{j}{(z-\lambda_j)_+^\rho}$$ for $\rho>1$ have been studied
for the Dirichlet problem. We obtain these means from
Corollary~\ref{corollary to generalized Riesz means} choosing
$f(\lambda)=(1-\lambda)^{\rho}$ and putting $z=1/t$. We therefore
have

\begin{corollary}\label{generalization of HaHe}
Suppose that the pair $(H,G)$ satisfy the relations
\eqref{commutator assumptions}.  Then
\begin{equation}\label{HHgen}
\frac{R_{\rho}(z)}{z^{\rho+2\gamma/\beta}}
\end{equation}
is a nondecreasing function for $0 < z < \inf{\sigma_{ess}(H)}$.
\end{corollary}

Formula \eqref{HHgen} is identical in form to an inequality in
\cite{HaHe}; the constant $2 \gamma/\beta$ has simply replaced
$d/2$ in the earlier article.  As a consequence we obtain a
Weyl-sharp bound,
\begin{equation}\label{Transformed bound}
\frac{\overline{\lambda_n}}{\overline{\lambda_k}}  \le
\left(\frac{2 (\beta+\gamma)^{1+\frac{\beta}{2 \gamma}}}{\beta
(\beta + 2 \gamma)^{\frac{\beta}{2 \gamma}}}\right)
\left(\frac{n}{k}\right)^{\frac{\beta}{2 \gamma}},
\end{equation}
provided that $n \ge \left(1 + \frac{2 \gamma}{\beta}\right) k$.
(In case the constant $\alpha$ has not been set to $0$, the
averages $\overline{\lambda_{n,k}}$ on the left side are both
replaced by $\overline{\lambda_{n,k}} + \alpha/\beta$.)

\section{Application to the Dirichlet Laplacian}
For the Dirichlet Laplacian $H=-\Delta_D$ on a domain $\Omega$ in
$\mathbb{R}^d$ such that $H$ has only eigenvalues, we put
$G=x_{\ell}$, the multiplication operator by a suitable Cartesian
coordinate. As shown for example in \cite{HaSt} we then have
$\alpha=0$, $\beta=\frac{4}{d}$ and $\gamma=1$. For the Dirichlet
Laplacian the inequality \eqref{T2} of our main 
Theorem~\ref{inequalities for bottom of spectrum} reads
\begin{equation}\label{T2-Dirichlet}
\sum_{j=1}^{n}\left(f(\lambda_j)+\frac{2}{d}\lambda_jf'(\lambda_j)\right)\leq
    n\left(f(\lambda_n)+\frac1{2}f'(\lambda_n)(\lambda_{n+1}-\lambda_n)\right).
\end{equation}
We claim that all estimates of the form \eqref{T2-Dirichlet} are
sharp in the semiclassical limit. Indeed, recall that according to
the Weyl law, on any bounded domain the semiclassical limit of the
eigenvalue $\lambda_n$ is given by
\begin{equation*}
   \lambda_n\sim
    C_d\bigg(\frac{n}{V}\bigg)^{\frac{2}{d}}
\end{equation*}
as $n\rightarrow\infty$, where $V$ denotes the volume of $\Omega$.
In terms of the counting function $N(\lambda)$ this is equivalent
to
\begin{equation}
    N(\lambda)\sim C_d^{-d/2}V\lambda^{d/2}.
\end{equation}
Now, for any function $f$, we have
\begin{align}
    \sum_{j=1}^{N(\lambda)}
    f(\lambda_j)&=N(\lambda)f(\lambda)-\int_0^{\lambda}f'(t)N(t)\;dt\\
    &\sim C_d^{-d/2}V
    \left(f(\lambda)-\int_0^{\lambda}f'(t)t^{d/2}\;dt\right),\nonumber
\end{align}
from which it easily follows that
\begin{align}
&\sum_{j=1}^{N(\lambda)}\left(f(\lambda_j)+\frac{2}{d}\lambda_jf'(\lambda_j)-f(\lambda_n)\right)\\
&\quad \sim C_d^{-d/2}V\left(\frac{2}{d}\lambda^{d/2+1}
f'(\lambda)-\int_0^{\lambda}
\left(1+\frac{2}{d}\right)f'(t)t^{d/2}+\frac{2}{d}f''(t)t^{d/2+1}
    \;dt\right)\nonumber\\
&\quad =\frac{2}{d}C_d^{-d/2}V\left(\lambda^{d/2+1}
f'(\lambda)-\int_0^{\lambda}(f'(t)t^{d/2+1})'\;dt\right)\nonumber\\
&\quad \sim o(N(\lambda)f(\lambda)).\nonumber
\end{align}
Consequently, Theorem~\ref{main theorem on moments} for the
moments
\begin{equation*}
     S_n(r)=\left(1+\frac{2
    r}{d}\right)\frac1{n}\sum_{j=1}^{n}\lambda_j^r
\end{equation*}
is sharp. Since the semiclassical limit of $S_n(r)^{\frac1{r}}$
does not depend on $r$ and is given by
\begin{equation*}
    S_n(r)^{\frac1{r}}\sim
    C_d\bigg(\frac{n}{V}\bigg)^{\frac{2}{d}}
\end{equation*}
as $n\rightarrow\infty$,
phase-space bounds on $S_n(r)$ for $r\leq 1$
follow from the Berezin-Li-Yau bound \cite{Be,LW1,LiYa} for $r=1$
by the monotonicity property \eqref{means small p}.

We 
can further
refine Theorem~\ref{main theorem on moments} for the
Dirichlet Laplacian by exploiting the right side of
inequality \eqref{T2} of Theorem~\ref{inequalities for bottom of
spectrum}, or respectively Corollary~\ref{corollary to inequalities
for bottom of spectrum},
in order to relate the arithmetic and geometric means of eigenvalues 
to the sizes of the eigenvalue gaps.  We begin by choosing
\begin{equation}\label{small powers and ln}
 f_z(\lambda)=\lambda^p-pz^{p}\ln \lambda + pz^{p}\ln z-z^{p}
\end{equation}
for $0<p\leq 3$, which corresponds to the choice $q=0$ in
\eqref{small powers}. Therefore $f_z(\lambda)$ in \eqref{small
powers and ln} satisfies the assumptions of 
Theorem~\ref{inequalities for bottom of spectrum} and 
Corollary~\ref{corollary to inequalities for bottom of spectrum} for all $z$
in an interval of the form $[\lambda_n,\Lambda_{n}(p)]$ for some
$\Lambda_{n}(p)$ which is determined by the condition
$f_z'(\lambda_n)+f_z''(\lambda_n)(\lambda_{n+1}-\lambda_n)\geq 0$.
Defining
\begin{equation*}
    \gamma_n := \frac{\lambda_{n+1}-\lambda_n}{\lambda_n}
\end{equation*}
we find
\begin{equation*}
    \Lambda_{n}(p)^p=\begin{cases}
\lambda_{n}^p\frac{1+(p-1)\gamma_n}{1-\gamma_n}&\text{if}\; \gamma_n<1\\ \infty&\text{otherwise.}\\
    \end{cases}
\end{equation*}
Defining
\begin{equation*}
    F_n(z) := npz^p\ln z-nz^p-npz^p\ln (e^{\frac{2}{d}}G_n)+nS_n(p)
\end{equation*}
and
\begin{equation*}
    \tilde{F}_n(z) := nf_z(\lambda_n)+n\frac1{2}f'_z(\lambda_n)(\lambda_{n+1}-\lambda_n),
\end{equation*}
we have
\begin{equation*}
    F_n(z)\leq \tilde{F}_n(z)
\end{equation*}
for all $z\in [\lambda_n,\Lambda_{n}(p)]$. We see that
$\tilde{F}_n(z)$ has a global minimum at
$\tilde{z}=\lambda_ne^{\gamma_n/2}\leq \Lambda_{n}(p)$. 
As the global minimum of $F_n(z)$ is below the global minimum of
$\tilde{F}_n(z)$, we obtain the inequality
\begin{equation}\label{p-ln-1}
    S_n(p)-(e^{\frac{2}{d}} G_n)^p\leq
    \lambda_n^p\left(1+\frac{p\gamma_n}{2}-e^{\frac{p\gamma_n}{2}}\right).
\end{equation}
We note that the left side of \eqref{p-ln-1} can be bounded
above by $-\frac1{2}(\frac{p\gamma_n}{2})^2$ which yields an
explicit upper bound on the gap $\lambda_{n+1}-\lambda_n$.
However, we find it more convenient to optimize the inequality
\begin{equation*}
    z^{-p}F_n(z)\leq z^{-p}\tilde{F}_n(z)
\end{equation*}
with respect to $z$. The right side has then a global minimum at
$\tilde{z}=\lambda_n(1+\frac{p\gamma_n}{2})^{\frac1{p}}\leq
\Lambda_{n}(p)$, while the left side has its global minimum at
$z=S_n(p)^{\frac1{p}}$. After taking the exponential on both
sides we therefore obtain the inequality
\begin{equation}\label{p-ln-2}
    S^{\frac1{p}}_n(p)\big/(e^{\frac{2}{d}} G_n)\leq
   \left(1+\frac{p\gamma_n}{2}\right)^{\frac1{p}}\big/e^{\frac{\gamma_n}{2}}.
\end{equation}
Extending the above discussion to all pairs $(p,q)$ of 
Theorem~\ref{main theorem on moments} we obtain the following refinement
for the Dirichlet Laplacian:

\begin{theorem}\label{main theorem on moments Dirichlet Laplacian} Let $p>q>0$ such that
$q\leq \min(1,p)$ and $p+q\leq 3$. Then for all $n$ we have
\begin{equation}\label{means small p Dirichlet Laplacian}
    S_n(p)^{\frac1{p}}\big/\left(1+\frac{p}{2}\gamma_n\right)^{\frac1{p}}
\leq  S_n(q)^{\frac1{q}}\big/\left(1+\frac{q}{2}\gamma_n\right)^{\frac1{q}}.
\end{equation}
In particular, for $0<p\leq 1$ the function
$$p\mapsto
S_n(p)^{\frac1{p}}\big/\left(1+\frac{p}{2}\gamma_n\right)^{\frac1{p}}$$ is
nonincreasing.
\end{theorem}

\begin{acknowledgements}  The hospitality of the \'Ecole Polytechnique F\'ed\'erale de Lausanne
for E.H. in support of this collaboration is gratefully acknowledged.
\end{acknowledgements}

\begin{thebibliography}{999}

\bibitem{As}
Ashbaugh, M.S.: The universal eigenvalue bounds of
Payne-P\'{o}lya-Weinberger, Hile-Protter, and H. C. Yang. In
{Spectral and inverse spectral theory (Goa, 2000)},
{Proc.\ Indian Acad.\ Sci.\ Math.\ Sci.} {112},
3--30 (2002)


\bibitem{AsHe}
Ashbaugh, M.S., Hermi, L.: A unified approach to universal
inequalities for eigenvalues of elliptic operators.
{Pacific J.  Math.} {217}, 201--219 (2004)


\bibitem{Be} Berezin, F.: Convariant and contravariant
symbols of operators. {Izv.  Akad. \ Nauk SSSR} {37},
1134--1167 (1972)
[In Russian.  English transl. in Math. \ USSR-Izv.
{6} (1972) 1117-1151 (1973)]

\bibitem{BeJa}  Bethe, H.A.  and Jackiw, R.W.:
{Intermediate quantum mechanics}, 2nd edn.
W.A. Benjamin, New York (1968)

\bibitem{ChYa}
Cheng, Q.M., Yang, H.C.: Estimates on eigenvalues of Laplacian.
{Math. Ann.} {331}, 445--460 (2005)


\bibitem{CY1} 
Cheng, Q.M. and Yang, H.C.: Estimates on eigenvalues of Laplacian.
{Math. Ann. 331}, 445--460 (2005)

\bibitem{CY2} Cheng, Q.M. and Yang, H.C.:
Inequalities for eigenvalues of Laplacian on domains and complex 
hypersurfaces in complex projective spaces. {J. Math. Soc. Japan 58},
545--561 (2006)


\bibitem{EHI}
El Soufi, A., Harrell, E.M. and Ilias, S.: Universal inequalities
for the eigenvalues of Laplace and Schr\"odinger operators on
submanifolds. 
{Trans. Amer. Math. Soc.}, to appear

\bibitem{Ha07}
Harrell, E.M.: Commutators, eigenvalue gaps, and mean curvature in
the theory of Schr\"odinger operators. {Commun.\ Part.\
Diff.\ Eq.} {32}, 401--413 (2007)


\bibitem{HaHe}
Harrell, E.M. and Hermi, L.: Differential inequalities for Riesz
means and Weyl-type bounds for eigenvalues.
{J. Funct. Analysis}
{254}, 3173--3191 (2008)


\bibitem{HaHe2}
Harrell II, E.M. and Hermi, L.: On Riesz Means of Eigenvalues.
preprint 2007


\bibitem{HaSt} Harrell II, E.M. and Stubbe, J.: On trace identities
and universal eigenvalue estimates for some partial differential
operators.  {Trans.\ Amer.\ Math.\ Soc.}  {349}, 1797--1809  (1997)
%


\bibitem{Lap} Laptev, A.: Dirichlet and Neumann
eigenvalue problems on domains in Euclidean spaces. {J.\
Funct.\ Anal.} {151}, 531--545 (1997)

\bibitem{LW1} 
Laptev, A. and Weidl, T.: Recent results on Lieb-Thirring
inequalities.  {Journ\'ees ``\'Equations aux D\'eriv\'ees Partielles'' 
(La Chapelle sur Erdre, 2000)},
Exp. No. XX, 14 pp., Univ.\ Nantes, Nantes, 2000.

\bibitem{LW2} Laptev, A. and Weidl, T.:
Sharp Lieb-Thirring inequalities in high dimensions. {Acta
Math.} {184}, 87--111 (2000)

\bibitem{LePa}
Levitin, M.  and Parnovski, L.: Commutators, spectral trace
identities, and universal estimates for eigenvalues. {J.\
Funct.\ Anal.} 192, 425--445 (2002)

\bibitem{Li}
Li, P.: Eigenvalue estimates on homogeneous manifolds.
{Comment.  Math. Helvitici} {55}, 347-363 (1980)


\bibitem{LiYa} Li, P. and Yau, S.-T.: On the Schr\"odinger equation and
the eigenvalue problem. 
{Commun.\ Math.\ Phys.} {88}, 309--318 (1983)

\bibitem{ReSi}  Reed, M. and Simon, B.:
{Methods of Modern Mathematical Physics, {I}:
Functional Analysis}. Revised and enlarged ed. Academic Press, New York (1980)

\bibitem{ReSiII}  Reed, M. and Simon, B.: 
{Methods of Modern Mathematical Physics, {II}: Fourier
Analysis, Self-Adjointness}. Academic Press, New York (1975)

\bibitem{We}
Webster, R.: {Convexity}.  Oxford University Press, Oxford (1994)


\bibitem{Ya}
Yang, H.C.: Estimates of the difference between
consecutive eigenvalues, 1995 preprint (revision of International
Centre for Theoretical Physics preprint IC/91/60, Trieste, Italy,
April 1991)


\end{thebibliography}

\end{document}